\documentclass[reqno,10pt]{article}
\usepackage{hyperref}
\usepackage{amsmath}
\usepackage{amssymb}
\usepackage{mathrsfs}

\newtheorem{thm}{Theorem}[section]
\newtheorem{lemma}[thm]{Lemma}
\newtheorem{corol}[thm]{Corollary}
\newtheorem{propos}[thm]{Proposition}
\newtheorem{rema}{Remark}[section]
\def\bp{\begin{propos}}
\def\ep{\end{propos}}
\def\bt{\begin{thm}}
\def\et{\end{thm}}
\def\bco{\begin{corol}}
\def\eco{\end{corol}}
\def\bl{\begin{lemma}}
\def\el{\end{lemma}}
\def\br{\begin{rema}}
\def\er{\end{rema}}
\def\be{\begin{equation}}
\def\ee{\end{equation}}
\def\ba{\begin{array}}
\def\ea{\end{array}}

\def\P{{\mathbb P}}
\def\E{{\mathbb E}}

\def\Z{{\mathbb Z}}

\def\QED{\hfill$\square$\vskip 3mm}

\def\Dp{\displaystyle}

\begin{document}

\title{\Large ON A LOWER BOUND FOR THE TIME CONSTANT OF\\[2mm] FIRST-PASSAGE PERCOLATION\\[6mm]
}
\author{{ Xian-Yuan Wu}$^1$\thanks{Research supported in part
by the Natural Science Foundation of China (No. 10301023) }\and
{Ping Feng}$^*$}
\date{}
\maketitle

\begin{abstract}
We consider the Bernoulli first-passage percolation on $\mathbb Z^d\
(d\ge 2)$. That is, the edge passage time is taken independently to
be 1 with probability $1-p$ and 0 otherwise. Let ${\mu(p)}$ be the
time constant. We prove in this paper that \[ \mu(p_1)-\mu({p_2})\ge
\frac{\mu(p_2)}{1-p_2}(p_2-p_1)\] for all $ 0\leq p_1<p_2< 1$ by
using Russo's formula.

\vskip 3mm


\end{abstract}

\section{Introduction and statement of the results.}
\renewcommand{\theequation}{1.\arabic{equation}}
\setcounter{equation}{0}

We begin with the general first-passage percolation on $\mathbb
Z^d$. Let $\{t(e): e\in\mathbb Z^d\}$ be a sequence of i.i.d.
positive random variables with common distribution $F$, $t(e)$ is
the random passage time of edge $e$ and $F$ is the edge-passage
distribution of the model. For any path
$\gamma=\{e_1,e_2,\ldots,e_n\}$, the passage time of $\gamma$ is
\[T(\gamma):=\sum_{t=1}^n t(e_k).\] For any vertices $u,v\in\mathbb Z^d$ and vertex sets $A,B\subset\mathbb Z^d$,
let
\[T(u,v):=\inf_{\gamma \ni u,v}T(\gamma);\ \ T(A,B):=\inf_{u\in A,v\in B}T(u,v)\]
be the passage time from $u$ to $v$ and the passage time from $A$ to
$B$.

Let $0$ be the origin of $\mathbb Z^d$, $\hat
e_1=(1,0,\ldots,0)\in\mathbb Z^d$ and
$H_n=\{u=(u_1,u_2,\ldots,u_n)\in\mathbb Z^d:u_1=n\}$. Define

\[a_{0,n}:=T(0,n\hat e_1),\ \ \ b_{0,n}:=T(0,H_n).\]
To restrict $a_{0,n},b_{0,n}$ on cylinders, let
$$\Gamma^{cyl}(0,n\hat e_1)=\{\gamma:0,\ n\hat e_1\in \gamma \mbox{
and }\forall\ u\in\gamma, 0\leq u_1\leq n\}\ \ \ \ \ \ \ \ \ \ \
$$
$$\Gamma^{cyl}(0,H_n)=\{\gamma:0\in\gamma,\ \gamma\cap
H_n\not=\emptyset,\mbox{ and }\forall\ u\in\gamma, 0\leq u_1\leq
n\}$$ and define $$ t_{0,n}:=\Dp\inf_{\gamma\in \Gamma^{cyl}(0,n\hat
e_1)}T(\gamma);\ \ s_{0,n}:=\Dp\inf_{\gamma\in
\Gamma^{cyl}(0,H_n)}T(\gamma).$$

The time constant $\mu$ of the model is the common limit of
$\theta_{0,n}/n$ when ${n\rightarrow\infty}$ for $\theta=a,b,t$ or
$s$. Here we will not introduce all the detailed situations for the
above convergence under various moment conditions of $F$, and only
point out that, in most cases, for $\theta=a,b,t$ or $s$,

\be\label{1.1}\frac {\theta_{0,n}}n\rightarrow\mu=\mu(F)\ \
\mbox{a.s.} \ \mbox{as} \ n\rightarrow\infty.\ee For the details on
the convergence to $\mu$, one may refer to \cite{1,1.3,1.2,1.4}.

It is straightforward that $\theta_{0,n}, \theta=a,b,t$ or $s$,
depends on the states of infinitely many edges. The following is
another limit representation of $\mu $ given by Grimmett and Kesten
\cite{2}, from which, $\mu$ is represented as the limit of random
variables which only depend on the states of finitely many edges.

For any fixed $n\ge 1$, let $B_n=\{u\in\mathbb Z^d:0\leq u_i\leq n,
1\leq i\leq d\}$ be the box with side length $n$. Let
\[\phi_{0,n}=\inf\{T(\gamma):\gamma \mbox{ is a path in }B_n \mbox{ from }\{0\}\times[0,n]^{d-1}\mbox{ to }\{n\}\times[0,n]^{d-1}\}\]

Grimmett and Kesten \cite{2} proved that, if the time-passage
distribution $F$ satisfying: \[ \int(1-F(x))^4dx<\infty \mbox{ for
}d=2;\mbox { or}\ \int x^2dF(x)<\infty\mbox{ for }d\geq 3\] then
\be\label{1.2}\frac {\phi_{0,n}}n\rightarrow\mu \ \ \mbox{a.s. and
in}\ L^1,\ \ \mbox{as} \ n\rightarrow\infty. \ee

The first problem for time constant $\mu$ is: when will $\mu>0$?
Kesten \cite{3} solved this problem for all $d\geq 2$ as:
\be\label{1.3} \mu>0
 \Leftrightarrow F(0)<p_c(d),\ee where $F(0)=\mathbb P(t(e)=0)$ and $p_c(d)$ be the
critical probability for the general bond percolation on $\mathbb
Z^d$.

Further study on $\mu$ is carried out to solve such a problem: How
does $\mu=\mu(F)$ depend on the edge-passage distribution $F$? Berg
and Kesten \cite{4} solved this problem in part. As our result is a
further research in this direction, in the next paragraph, we
introduce the results of Berg and Kesten in detail.

Let's begin with some notations. For any given edge-passage
distributions $F$, let ${\rm supp}(F)=\{x\geq 0:F(x)>0\}$ be the
support of $F$, let $\lambda(F)=\inf{\rm supp}(F)$. We say $F$ is
{\it useful}, if
\[\lambda(F)=0 \mbox{ and } F(0)<p_c(d),  \mbox{ or }\lambda(F)>0 \mbox{
and }F(\lambda)<\vec{p}_c(d), \]where $\vec{p}_c(d)$ is the critical
probability for directed bond percolation on $\mathbb Z^d$. For two
edge-passage distributions $F$ and $\tilde F$, we say $\tilde F$ is
more {\it variable} than $F$, if \be\label{1.4}\int\varphi(x)d\tilde
F(x)\leq \int\varphi(x)dF(x)\ee for all increasing convex function
$\varphi$. Clearly, by the above definition, ``$\tilde F$ is more
{\it variable} than $F$" is a weaker condition than ``$\tilde F$ is
{\it stochastically dominated} by $F$", note that the latter
requires equation (\ref{1.4}) hold for all increasing $\varphi$.

\bt\label{a}${\mbox{\rm [Berg and Kesten \cite{4}]}}$\\ (a) Let $F$
and $\tilde F$ be two edge-passage distribution functions, if
$\tilde F$ is more {\it variable} than $F$, then
\[\mu(\tilde F)\leq \mu(F);\]
(b) if, in addition, $F$ is useful and $F\not=\tilde F$, then
\[\mu(\tilde F)<\mu(F).\]\et

Theorem \ref{a} gives sufficient conditions for (strict) inequality
between $\mu(\tilde F)$ and $\mu(F)$, but for the difference
$\mu(F)-\mu(\tilde F)$, no information is provided. One may ask:
what can we say for such a difference? In this paper, for the
simplest case, i.e., under the following Bernoulli setting, we give
a nontrivial lower bound for this difference.

From now on, we take $\{t(e):e\in\Z^d\}$ to be the i.i.d. random
variable sequence such that $t(e)=1$ with probability $1-p$ and
$t(e)=0$ with probability $p$, $p\in [0,1]$. Write $\P_p$ as the
percolation measure and $\E_p$ as its expectation. Write $\mu(p)$ as
the corresponding time constant. By (\ref{1.3}) and Theorem \ref{a},
$\mu(p)$ decreases strictly in $p$ when $p\in [0,p_c(d))$, i.e.,
\be\label{1.5}\frac{\mu(p_1)}{\mu(p_2)}>1\ee for all $0\leq
p_1<p_2<p_c(d)$.

\vskip 5mm Now, we state our main result as follows.

\bt\label{m} For the above Bernoulli first-passage percolation
model, let $\mu(p)$ be its time constant. We have that
${\mu(p)}/{(1-p)}$ decreases in $p$ and then
\be\label{1.6}\mu(p_1)-\mu({p_2})\ge\frac{\mu(p_2)}{1-p_2}(p_2-p_1)\ee
for all $0\le p_1<p_2\le1$.\et

\br By the monotonicity of ${\mu(p)}/{(1-p)}$ and (\ref{1.3}), when
$0\leq p_1<p_2<p_c(d)$, one has
\be\label{1.7}\frac{\mu(p_1)}{\mu(p_2)}\ge 1+\frac{p_2-p_1}{1-p_2}.
\ee This is a concretion of (\ref{1.5}).\er


\section{Proof of Theorem \ref{m}}
\renewcommand{\theequation}{2.\arabic{equation}}
\setcounter{equation}{0}

To use the Russo's formula, we first give the definition of {\it
pivotal} edges according to Grimmett \cite{5}. For any edge $e$ and
configuration $\omega$, let $\omega_e$ be the configuration such
that $\omega_e(f)=\omega(f)$ for all $f\not= e$ and
$\omega_e(e)=1-\omega(e)$.

Recall that $B_n=[0,n]^d\cap\Z^d$. Suppose that $A$ be an event
which only depends on edges of $B_n$. We say edge $e\in B_n$ is {\it
pivotal} for pair $(A,\omega)$, if
\[I_A(\omega)\not=I_A(\omega_e),\] where $I_A$ be the indicator function of $A$. Write $S_e(A)$ as the event
that $e$ is a {pivotal } edge for $A$, i.e. \be
S_e(A)=\{\omega:e\mbox{ is pivotal for pair }(A,\omega)\}.\ee By the
above definition, $S_e(A)$ is independent of $t(e)$. Denote by
$N(A)$ the number of {pivotal} edges of $A$, i.e. \be
N(A)(\omega)=|\{e\in B_n:\omega\in S_e(A)\}|.\ee Event $A$ is called
{\it increasing} if $\omega\in A$ and $\omega\leq\omega'$ imply
$\omega'\in A$, where $\omega\leq\omega'$ means $\omega(e)\leq
\omega'(e)$ for all $e$. The Russo's formula says that (in our
setting), if $A$ is increasing, then \be
\frac{d\P_p(A)}{dp}=-\E_p(N(A)).\ee

\vskip 5mm\noindent{\it Proof of Theorem \ref{m}}: Firstly, by
equation (\ref{1.2}), we have \be\label{2.1}\mu(p)=
\lim\limits_{n\rightarrow\infty}\frac{E_{p}\phi_{0,n}}{n}\ee for all
$p\in[0,1]$.

For any integer $k\geq 1$, let $A_{n,k}=\{\phi_{0,n}\geq k\}$.
Clearly, $A_{n,k}$ is increasing and only depends on edges in $B_n$.
Rewrite $\E_p(\phi_{0,n})$ as
\be\label{2.2}\E_p(\phi_{0,n})=\sum_{k=1}^{\infty}\P_p(A_{n,k}).\ee

For any $0\leq p_1<p_2\leq 1$, by (\ref{2.1}) and (\ref{2.2}), we
have \be\label{2.3}
\ba{rl}\mu(p_1)-\mu(p_2)&=\Dp\lim_{n\rightarrow\infty}\frac 1n
\sum_{k=1}^\infty(\P_{p_1}(A_{n,k})-\P_{p_2}(A_{n,k}))\\[5mm]
&=\Dp\lim_{n\rightarrow\infty}\frac 1n
\sum_{k=1}^\infty\int_{p_1}^{p_2}-\frac {d\P_p(A_{n,k})}{dp}dp.\\
\ea \ee Using the Russo's formula and the fact that $A_{n,k}$ is
increasing, we have \be\label{2.4}
\ba{ll}&\displaystyle\frac {d\P_{p}(A_{n,k})}{dp}=-\E_{p}(N(A_{n,k}))=-\Dp\sum_{e\in B_n}\P_{p}(S_e(A_{n,k}))\\[2mm]
&=-\Dp\frac{1}{1-p}\sum_{e\in B_n}\P_{p}(\{t(e)=1\}\cap S_e(A_{n,k}))\\[2mm]
&=-\Dp\frac{1}{1-p}\sum_{e\in B_n}\P_{p}(A_{n,k}\cap S_e(A_{n,k}))\\[2mm]
&=-\Dp\frac{1}{1-p}\sum_{e\in B_n}\P_{p}(S_e(A_{n,k})\mid A_{n,k})\P_{p}(A_{n,k})\\[2mm]
&=-\Dp\frac{1}{1-p}\E_{p}(N(A_{n,k})\mid A_{n,k})\P_{p}(A_{n,k}).\\
\ea \ee Note that the third equality comes from the independence of
$t(e)$ and $S_e(A_{n,k})$.

To finish the proof, we have to give appropriate lower bound for
$\E_{p}(N(A_{n,k})\mid A_{n,k})$. To this end, for any configuration
$\omega\in A_{n,k}$, we give lower bounds to $N(A_{n,k})(\omega)$ in
the following two cases respectively: 1) $\phi_{0,n}(\omega)\geq
k+1$; 2) $\phi_{0,n}(\omega)=k$.

We first deal with the case of $\phi_{0,n}(\omega)\geq k+1$. For any
$e\in B_n$, because $\omega_e$ only differs from $\omega$ in edge
$e$, the change from $\omega$ to $\omega_e$ can at most decrease
$\phi_{0,n}$ by 1, this implies that $\phi_{0,n}(\omega_e)\geq k$,
and $\omega_e\in A_{n,k}$. By the definition of pivotal edges, we
know that $e$ is not pivotal for $(A_{n,k},\omega)$. So
\be\label{2.5}N(A_{n,k})(\omega)=0.\ee

Now, we consider the case of $\phi_{0,n}(\omega)= k$. For any $e\in
B_n$, if $e$ is pivotal for $(A_{n,k},\omega)$, we declare that
$\omega(e)=1$. Actually, if $\omega(e)=0$, then the change from
$\omega$ to $\omega_e$ will increase $\phi_{0,n}$, so we have
$\phi_{0,n}(\omega_e)\geq\phi_{0,n}(\omega)=k$ and $\omega_e\in
A_{n,k}$, this leads to a contradiction.

Suppose $\gamma$ be a path in $B_n$ from $\{0\}\times [0,n]^{d-1}$
to $\{n\}\times [0,n]^{d-1}$ with
$T(\gamma)(\omega)=\phi_{0,n}(\omega)=k$. If $e\in\gamma$ satisfying
$\omega(e)=1$, then $T(\gamma)(\omega_e)=k-1$. This implies that
$\phi_{0,n}(\omega_e)\leq k-1$ and $\omega_e\notin A_{n,k}$. Thus,
by the definition of pivotal edges, $e$ is pivotal for pair
$(A_{n,k},\omega)$.

By the arguments in the last two paragraphs, we have
\be\label{2.6}N(A_{n,k})(\omega)\geq
|\{e\in\gamma:\omega(e)=1\}|=T(\gamma)(\omega)=k\ee for all
$\omega\in A_{n,k}$.

Combining (\ref{2.5}) and (\ref{2.6}), we have
\be\label{2.7}\begin{array}{ll} \E_{p}(N(A_{n,k})\mid
A_{n,k})\P_{p}(A_{n,k})&=\displaystyle\sum_{\omega\in
A_{n,k}}N(A_{n,k})(\omega)\cdot
\frac{\P_{p}(\omega)}{\P_{p}(A_{n,k})}\P_{p}(A_{n,k})\\[5mm]
&\geq\displaystyle\sum_{\{\omega:\phi_{0,n}(\omega)=k\}}k\cdot \P_{p}(\omega)\\[5mm]
&=k\cdot \P_{p}(\{\omega:\phi_{0,n}(\omega)=k\}).
\end{array}\ee
Finally, by (\ref{2.3}), (\ref{2.4}) and (\ref{2.7}), using the
Fubini's theorem and the Fatou's lemma, we have \be\label{2.8}
\begin{array}{ll}
\mu(p_1)-\mu(p_2)&=\displaystyle\lim_{n\rightarrow\infty}\frac 1n\sum_{k=1}^\infty\int_{p_1}^{p_2}\frac 1{1-p}\E_{p}(N(A_{n,k})\mid A_{n,k})\P_{p}(A_{n,k})dp\\[5mm]
&\geq\displaystyle\lim_{n\rightarrow\infty}\frac 1n\sum_{k=1}^\infty\int_{p_1}^{p_2}\frac 1{1-p}k\cdot \P_{p}(\{\omega:\phi_{0,n}(\omega)=k\})dp\\[5mm]
&=\displaystyle\lim_{n\rightarrow\infty}\frac 1n\int_{p_1}^{p_2}\frac 1{1-p}\sum_{k=1}^\infty k\cdot \P_{p}(\{\omega:\phi_{0,n}(\omega)=k\})dp\\[5mm]
&=\displaystyle\lim_{n\rightarrow\infty}\frac 1n\int_{p_1}^{p_2}\frac 1{1-p}\E_{p}(\phi_{0,n})dp\\[5mm]
&\geq\displaystyle\int_{p_1}^{p_2}\frac 1{1-p}\liminf_{n\rightarrow\infty}\frac 1n\E_{p}(\phi_{0,n})dp\\[5mm]
&=\displaystyle\int_{p_1}^{p_2}\frac {\mu(p)} {1-p}\,dp\\
\end{array}
\ee for all $0\leq p_1<p_2<1$. Clearly, the inequality (\ref{2.8})
is equivalent to the following differential inequality
\be\label{2.9}\frac{d[{\mu(p)}/{(1-p})]}{dp}\leq 0,\ \ 0\leq p<1.\ee
This gives that \[ \int_{p_1}^{p_2}\frac
{\mu(p)}{1-p}\,dp\ge\frac{\mu(p_2)}{1-p_2}(p_2-p_1)\] for all $0\leq
p_1<p_2<1$ and we finish the proof of Theorem \ref{m}. \QED

\bigskip

\noindent{\bf Acknowledgements.} This work was begun when one of us
(Xian-Yuan~Wu) was visiting Institute of Mathematics, Academia
Sinica. He is thankful to the probability group of IM-AS for
hospitality. We much thank Yu Zhang for drawing our attention to
this problem and for useful advices.

\vskip 5mm
\begin{minipage}{7cm}
{\noindent\hskip-3mm $^{~1}$Department of Mathematics, Capital
Normal University, 100037, Beijing, P. R. China.\\
E-mail: \texttt{wuxy@mail.cnu.edu.cn} }
\end{minipage}
\end{document}